
\def\LaTeX{%
  \let\Begin\begin
  \let\End\end
  \let\salta\relax
  \let\finqui\relax
  \let\futuro\relax}

\def\UK{\def\our{our}\let\sz s}
\def\USA{\def\our{or}\let\sz z}

\UK



\LaTeX

\USA


\salta

\documentclass[twoside,12pt]{article}
\setlength{\textheight}{24cm}
\setlength{\textwidth}{16cm}
\setlength{\oddsidemargin}{2mm}
\setlength{\evensidemargin}{2mm}
\setlength{\topmargin}{-15mm}
\parskip2mm


\usepackage[usenames,dvipsnames]{color}
\usepackage{amsmath}
\usepackage{amsthm}
\usepackage{amssymb}
\usepackage[mathcal]{euscript}
\usepackage{cite}

%
%


\definecolor{viola}{rgb}{0.3,0,0.7}
\definecolor{ciclamino}{rgb}{0.5,0,0.5}

\def\gianni #1{{\color{green}#1}}
\def\pier #1{{\color{red}#1}}
\def\juerg #1{{\color{viola}#1}}

\def\pier #1{#1}
\def\juerg #1{#1}
\def\gianni #1{#1}




\bibliographystyle{plain}


%

\finqui

\def\Beq{\Begin{equation}}
\def\Eeq{\End{equation}}
\def\Bsist{\Begin{eqnarray}}
\def\Esist{\End{eqnarray}}

\def\Bthm{\Begin{theorem}}
\def\Ethm{\End{theorem}}

\def\Bcor{\Begin{corollary}}
\def\Ecor{\End{corollary}}
\def\Brem{\Begin{remark}\rm}
\def\Erem{\End{remark}}

\def\Bcenter{\Begin{center}}
\def\Ecenter{\End{center}}
\let\non\nonumber




\def\step #1 \par{\medskip\noindent{\bf #1.}\quad}


\def\Lip{Lip\-schitz}
\def\Holder{H\"older}

\def\aand{\quad\hbox{and}\quad}

\def\lhs{left-hand side}
\def\rhs{right-hand side}

\def\omegalimit{$\omega$-limit}


\def\bhv{behavi\our}


\def\multibold #1{\def\arg{#1}%
  \ifx\arg\pto \let\next\relax
  \else
  \def\next{\expandafter
    \def\csname #1#1#1\endcsname{{\bf #1}}%
    \multibold}%
  \fi \next}

\def\pto{.}

\def\multical #1{\def\arg{#1}%
  \ifx\arg\pto \let\next\relax
  \else
  \def\next{\expandafter
    \def\csname cal#1\endcsname{{\cal #1}}%
    \multical}%
  \fi \next}


\def\multimathop #1 {\def\arg{#1}%
  \ifx\arg\pto \let\next\relax
  \else
  \def\next{\expandafter
    \def\csname #1\endcsname{\mathop{\rm #1}\nolimits}%
    \multimathop}%
  \fi \next}

\multibold
qwertyuiopasdfghjklzxcvbnmQWERTYUIOPASDFGHJKLZXCVBNM.

\multical
QWERTYUIOPASDFGHJKLZXCVBNM.

\multimathop
diag dist div dom mean meas sign supp .


\def\accorpa #1#2{\eqref{#1}--\eqref{#2}}
\def\Accorpa #1#2 #3 {\gdef #1{\eqref{#2}--\eqref{#3}}%
  \wlog{}\wlog{\string #1 -> #2 - #3}\wlog{}}


\def\separa{\noalign{\allowbreak}}

\def\graffe #1{\mathopen\{#1\mathclose\}}

\def\<#1>{\mathopen\langle #1\mathclose\rangle}
\def\norma #1{\mathopen \| #1\mathclose \|}

\def\[#1]{\mathopen\langle\!\langle #1\mathclose\rangle\!\rangle}

\def\iot {\int_0^t}

\def\intQt{\int_{Q_t}}

\def\iO{\int_\Omega}
\def\iG{\int_\Gamma}

\def\intSt{\int_{\Sigma_t}}
\def\intQi{\int_{\Qi}}
\def\intSi{\int_{\Si}}

\def\dt{\partial_t}
\def\dn{\partial_\nu}

\def\cpto{\,\cdot\,}

\def\checkmmode #1{\relax\ifmmode\hbox{#1}\else{#1}\fi}
\def\aeO{\checkmmode{a.e.\ in~$\Omega$}}
\def\aeQ{\checkmmode{a.e.\ in~$\QT$}}
\def\aeG{\checkmmode{a.e.\ on~$\Gamma$}}
\def\aeS{\checkmmode{a.e.\ on~$\ST$}}
\def\aet{\checkmmode{a.e.\ in~$(0,T)$}}

\def\aat{\checkmmode{for a.a.~$t\in(0,T)$}}

\def\Aet{\checkmmode{a.e.\ in~$(0,+\infty)$}}
\def\Aat{\checkmmode{for a.a.~$t\in(0,+\infty)$}}

\def\limn{\lim_{n\to\infty}}


\def\erre{{\mathbb{R}}}

\def\enne{{\mathbb{N}}}




\def\genspazio #1#2#3#4#5{#1^{#2}(#5,#4;#3)}
\def\spazio #1#2#3{\genspazio {#1}{#2}{#3}T0}

\def\L {\spazio L}
\def\H {\spazio H}
\def\W {\spazio W}

\def\C #1#2{C^{#1}([0,T];#2)}
\def\spazioinf #1#2#3{\genspazio {#1}{#2}{#3}{+\infty}0}
\def\LL {\spazioinf L}
\def\HH {\spazioinf H}

\def\CC #1#2{C^{#1}([0,+\infty);#2)}


\def\Lx #1{L^{#1}(\Omega)}
\def\Hx #1{H^{#1}(\Omega)}

\def\LxG #1{L^{#1}(\Gamma)}
\def\HxG #1{H^{#1}(\Gamma)}

\def\Luno{\Lx 1}
\def\Ldue{\Lx 2}

\def\Huno{\Hx 1}
\def\Hdue{\Hx 2}
\def\Hunoz{{H^1_0(\Omega)}}
\def\HunoG{\HxG 1}
\def\HdueG{\HxG 2}

\def\LunoG{\LxG 1}
\def\LdueG{\LxG 2}



\let\theta\vartheta
\let\eps\varepsilon
\let\phi\varphi

\let\hat\widehat

\let\TeXchi\chi                         
\newbox\chibox
\setbox0 \hbox{\mathsurround0pt $\TeXchi$}
\setbox\chibox \hbox{\raise\dp0 \box 0 }
\def\chi{\copy\chibox}



\def\suG{_{|\Gamma}}

\def\VG{V_\Gamma}
\def\HG{H_\Gamma}
\def\WG{W_\Gamma}
\def\nablaG{\nabla_\Gamma}
\def\DeltaG{\Delta_\Gamma}
\def\muG{\mu_\Gamma}
\def\rhoG{\rho_\Gamma}
\def\tauO{\tau_\Omega}
\def\tauG{\tau_\Gamma}
\def\fG{f_\Gamma}
\def\vG{v_\Gamma}
\def\wG{w_\Gamma}
\def\gG{g_\Gamma}
\def\zetaG{\zeta_\Gamma}
\def\xiG{\xi_\Gamma}
\def\jsoluz{((\mu,\muG),(\rho,\rhoG),(\zeta,\zetaG))}
\def\soluz{(\mu,\muG,\rho,\rhoG,\zeta,\zetaG)}
\def\soluzn{(\mun,\muGn,\rhon,\rhoGn,\zetan,\zetaGn)}
\def\soluzi{(\mui,\muGi,\rhoi,\rhoGi,\zetai,\zetaGi)}
\def\soluzs{(\rhos,\rhoGs,\zetas,\zetaGs)}
\def\Mu{(\mu,\muG)}
\def\Rho{(\rho,\rhoG)}
\def\Zeta{(\zeta,\zetaG)}
\def\Xi{(\xi,\xiG)}

\def\gstar{g^*}

\def\rhoz{\rho_0}
\def\rhoGz{{\rhoz}\suG}
\def\Rhoz{(\rhoz,\rhoGz)}
\def\mz{m_0}

\def\Beta{\hat\beta}
\def\BetaG{\Beta_\Gamma}
\def\betaG{\beta_\Gamma}
\def\betaeps{\beta_\eps}
\def\betaGeps{\beta_{\Gamma\!,\,\eps}}

\def\betaz{\beta^\circ}
\def\betaGz{\betaG^\circ}
\def\Pi{\hat\pi}
\def\PiG{\Pi_\Gamma}
\def\piG{\pi_\Gamma}

\def\mun{\mu^n}
\def\muGn{\mu_\Gamma^n}
\def\rhon{\rho^n}
\def\rhoGn{\rho_\Gamma^n}
\def\zetan{\zeta^n}
\def\zetaGn{\zeta_\Gamma^n}
\def\Mun{(\mun,\muGn)}
\def\Rhon{(\rhon,\rhoGn)}
\def\Zetan{(\zetan,\zetaGn)}
\def\un{u^n}

\def\rhoo{\rho^\omega}
\def\rhoGo{\rho_\Gamma^\omega}
\def\Rhoo{(\rhoo,\rhoGo)}

\def\mus{\mu^s}
\def\rhos{\rho^s}
\def\rhoGs{\rho_\Gamma^s}
\def\zetas{\zeta^s}
\def\zetaGs{\zeta_\Gamma^s}
\def\Rhos{(\rhos,\rhoGs)}
\def\Zetas{(\zeta,\zetaGs)}

\def\mui{\mu^\infty}
\def\muGi{\mu_\Gamma^\infty}
\def\rhoi{\rho^\infty}
\def\rhoGi{\rho_\Gamma^\infty}
\def\zetai{\zeta^\infty}
\def\zetaGi{\zeta_\Gamma^\infty}
\def\Mui{(\mui,\muGi)}
\def\Rhoi{(\rhoi,\rhoGi)}
\def\Zetai{(\zetai,\zetaGi)}

\def\musi{\mu_\infty}

\def\calVz{\calV_0}
\def\calHz{\calH_0}
\def\calVsz{\calV_{*0}}
\def\calVzp{\calV_0^{\,*}}
\def\calVp{\calV^{\,*}}

\def\calNO{\calN_\Omega}
\def\calNG{\calN_\Gamma}

\def\normaV #1{\norma{#1}_V}
\def\normaH #1{\norma{#1}_H}
\def\normaW #1{\norma{#1}_W}

\def\normaHG #1{\norma{#1}_{\HG}}

\def\normaHH #1{\norma{#1}_{\calH}}

\def\normaWW #1{\norma{#1}_{\calW}}

\def\Qi{Q_\infty}
\def\Si{\Sigma_\infty}
\def\QT{Q_T}
\def\ST{\Sigma_T}
\def\tn{t_n}

\def\longtime{longtime}
\def\Longtime{Longtime}
\def\CO{C_\Omega}
\def\CT{c_T}

\Begin{document}


%
\title{On the \longtime\ \bhv\ of a viscous \\[0.3cm] 
  Cahn--Hilliard system with convection\\[3mm] and dynamic boundary conditions}
\author{}
\date{}
\maketitle
\Bcenter
\vskip-1cm
{\large\sc Pierluigi Colli$^{(1)}$}\\
{\normalsize e-mail: {\tt pierluigi.colli@unipv.it}}\\[.25cm]
{\large\sc Gianni Gilardi$^{(1)}$}\\
{\normalsize e-mail: {\tt gianni.gilardi@unipv.it}}\\[.25cm]
{\large\sc J\"urgen Sprekels$^{(2)}$}\\
{\normalsize e-mail: {\tt sprekels@wias-berlin.de}}\\[.45cm]
$^{(1)}$
{\small Dipartimento di Matematica ``F. Casorati'', Universit\`a di Pavia}\\
{\small and Research Associate at the IMATI -- C.N.R. Pavia}\\
{\small via Ferrata 5, 27100 Pavia, Italy}\\[.2cm]
$^{(2)}$
{\small Department of Mathematics}\\
{\small Humboldt-Universit\"at zu Berlin}\\
{\small Unter den Linden 6, 10099 Berlin, Germany}\\[2mm]
{\small and}\\[2mm]
{\small Weierstrass Institute for Applied Analysis and Stochastics}\\
{\small Mohrenstrasse 39, 10117 Berlin, Germany}\\
[1cm]
{\it Dedicated to our friend Prof. Dr. Alexander Mielke\\[.1cm]
on the occasion of his 60th birthday
with best wishes}
\Ecenter
\Begin{abstract}
\juerg{In this paper, we study the longtime asymptotic behavior of a phase separation process
occurring in a three-dimensional domain containing a fluid flow of given velocity.
This process is modeled by a viscous convective Cahn--Hilliard system, which consists of 
two nonlinearly coupled second-order partial differential equations for the unknown
quantities, the chemical potential and an order parameter representing the scaled 
density of one of the phases. In contrast to other contributions, in which zero Neumann
boundary conditions were are assumed for both the chemical potential and the
order parameter, we consider the case of dynamic boundary conditions, which
model the situation when another phase transition takes place on the boundary.  The
phase transition processes in the bulk and on the boundary are driven by 
free energies functionals that may be nondifferentiable and have derivatives only in the
sense of (possibly set-valued) subdifferentials.  
For the resulting initial-boundary value system of Cahn--Hilliard type, general
well-posedness results have been established in \pier{a recent contribution by the same authors}. In the
present paper, we \pier{investigate} the asymptotic behavior of the solutions as times approaches infinity.      
More precisely, we study the \omegalimit\ (in~a suitable topology)
of~every solution trajectory.
Under the assumptions that the viscosity coefficients are strictly positive and that at least one of the underlying free energies is differentiable,
we prove that the \omegalimit\ is meaningful and that all of  its  elements
are  solutions to the corresponding stationary system,
where  the component representing the chemical potential is a constant.}

\End{abstract}
\salta
\pagestyle{myheadings}
\newcommand\testopari{\sc Colli \ --- \ Gilardi \ --- \ Sprekels}
\newcommand\testodispari{\sc Longtime behavior of a viscous convective Cahn--Hilliard system}
\markboth{\testodispari}{\testopari}
\finqui
%

\section{Introduction}
\label{Intro}
\setcounter{equation}{0}

\pier{The} recent paper \cite{CGS13} \pier{addresses}
an initial-boundary value problem for the Cahn--Hilliard system with convection
\Beq
  \dt\rho + \nabla\rho \cdot u - \Delta\mu = 0
  \aand
  \tauO \dt\rho - \Delta\rho + f'(\rho) = \mu
  \quad \hbox{in $\QT:=\Omega\times(0,T)$},
  \label{Isystem}
\Eeq
in the unknowns~$\rho$, the order parameter, and~$\mu$, the chemical potential.
In the above equations, 
$\tauO$~is a nonnegative constant,
$f'$~is the derivative of a double-well potential~$f$,
and $u$ is a given velocity field.

Typical and significant examples of $f$ 
are the so-called {\em classical regular potential}, the {\em logarithmic double-well potential\/},
and the {\em double obstacle potential\/}, which are given~by
\Bsist
  && f_{reg}(r) := \frac 14 \, (r^2-1)^2 \,,
  \quad r \in \erre, 
  \label{regpot}
  \\
  && f_{log}(r) := \bigl( (1+r)\ln (1+r)+(1-r)\ln (1-r) \bigr) - c_1 r^2 \,,
  \quad r \in (-1,1),
  \label{logpot}
  \\[2mm]
  && f_{2obs}(r) :=  c_2 \pier{(1 - r^2)} 
  \quad \hbox{if $|r|\leq1$}
  \aand
  f_{2obs}(r) := +\infty
  \quad \hbox{if $|r|>1$},
  \label{obspot}
\Esist
where the \pier{constants} in \eqref{logpot} and \eqref{obspot} satisfy
$c_1>1$ and $c_2>0$, 
so that that $f_{log}$ and $f_{2obs}$ are nonconvex.
In cases like \eqref{obspot}, one has to split $f$ into a nondifferentiable convex part 
(the~indicator function of $[-1,1]$ in the present example) and a smooth perturbation.
Accordingly, one has to replace the derivative of the convex part
by the subdifferential and interpret the second identity in \eqref{Isystem} as a differential inclusion.

As far as the conditions on the boundary $\Gamma:=\partial\Omega$ are concerned, 
instead of the classical homogeneous Neumann boundary conditions, 
the dynamic boundary condition for both $\mu$ and~$\rho$ are considered, namely,
\begin{align}
  &\pier{\dt\rhoG + \dn\mu - \DeltaG\muG = 0 \quad \hbox{and}}
  \non \\
  &\qquad\tauG \dt\rhoG + \dn\rho - \DeltaG\rhoG + \fG'(\rhoG) = \muG
 \quad \hbox{on $\ST:=\Gamma\times(0,T)$},
  \label{IdynBC}
\end{align}
where $\juerg{\muG=\mu_{|\Sigma_T}}$, $\juerg{\rhoG=\rho_{|\Sigma_T}}$, are the traces of $\mu$ and~$\rho$, respectively,
$\dn$~and $\DeltaG$ denote the outward normal derivative
and the Laplace--Beltrami operator on~$\Gamma$,
$\tauG$~is a nonnegative constant,
and $\fG'$ is the derivative of another potential~$\fG$.

\pier{The equations in \eqref{Isystem} aim to describe a class of evolution phenomena with phase separation
and fluid convection, in which the convection is represented by the term  $\nabla\rho \cdot u$ for a known  
velocity vector $u$. Regarding this system of partial differential equations, some boundary conditions 
are usually prescribed and the standard approach leads to the no-flux conditions
\Beq
  \dn\mu = 0, \quad \dn\rho = 0 
  \quad \hbox{ on $\Sigma_T$},
  \non
\Eeq
for both $\mu$ and $\rho$. On the contrary, as already announced, here we are interested to handle 
the dynamic boundary conditions \eqref{IdynBC}, which
set a  Cahn--Hilliard type system also on the boundary. The two potentials $f= \Beta +\Pi$ in the bulk and 
$\fG= \BetaG + \PiG$ on the boundary are both the sum of a convex and lower 
semicontinuous part and a (possibly concave) perturbation; they are not completely independent
but related by a suitable growth condition. Within the framework given by  \eqref{Isystem} and  
\eqref{IdynBC}, initial conditions should be prescribed for $\rho$ both in the bulk and on the boundary.}

\pier{All in all, the resulting initial and boundary value problem~reads
\begin{align}
  & \dt\rho + \nabla\rho \cdot u - \Delta\mu = 0
  \quad \hbox{in $Q_T$},
  \label{pier2}
    \\
  & \tauO \dt\rho - \Delta\rho + \beta(\rho) + \pi(\rho) \ni \mu
  \quad \hbox{in $Q_T$},
  \label{pier3}
  \\
  &\rhoG= \rho_{|_{\Sigma_T}}, \quad  \muG= \mu_{|_{\Sigma_T}} \aand  \dt\rhoG + \dn\mu - \DeltaG\muG = 0
  \quad \hbox{on $\Sigma_T$},
  \label{pier4}
  \\
  & \tauG \dt\rhoG + \dn\rho - \DeltaG\rhoG + \betaG(\rhoG) + \piG(\rhoG) \ni \muG
  \quad \hbox{on $\Sigma_T$},
  \label{pier5}
  \\
 & \rho(0) = \rhoz
  \quad \hbox{in $\Omega$} \aand \rhoG (0) = \rhoGz \quad \hbox{on $\Gamma$}.
    \label{pier6}
\end{align}}%
The paper \cite{CGS13} is devoted to the study
of the initial-boundary value problem \pier{\eqref{pier2}--\eqref{pier6}.} 
Under suitable assumptions and compatibility conditions on the potentials,
well-posedness and regularity results are proved.

\pier{The aim of the present paper is investigating} the \longtime\ \bhv.
More precisely, we study the \omegalimit\ (in~a suitable topology)
of~every trajectory~$\Rho$.
Under the additional assumptions that the viscosity coefficients $\tauO$ and $\tauG$
are strictly positive and that at least one of the potentials $f$ and $\fG$ is differentiable,
we prove that the \omegalimit\ is meaningful and that every element $\Rhoo$ of it
is a stationary solution $\Rhos$ of the system for $\Rho$
with some constant value $\mus$ of the chemical potential.

\pier{Let us now review some related literature. 
About Cahn--Hilliard problems, we quote the pioneering contributions  \cite{CahH, EllSh, NovCoh, bai, EllSt}
and observe that for this class of evolution processes it turns out that the phases do not diffuse\juerg{,} but  
they separately concentrate and form the so-called  spinodal decomposition.
A discussion on the modeling aspects of phase separation, spinodal decomposition and mobility 
of atoms between cells can be found in \cite{FG, Gu, Podio, CGPS3, CMZ11}). 
Up to our knowledge, in the case of a pure Cahn--Hilliard system, 
that is, with $\tauO=\tauG=0$, and without convective term ($u=0$), 
the problem \eqref{pier2}--\eqref{pier6} has been introduced
by Gal~\cite{G1} and formulated by {G}oldstein, {M}iranville and {S}chimperna~\cite{GMS11}. 
It has been studied from various viewpoints in other contributions 
(see~\cite{CGM13, CMZ11, CP14, GW, GMS11, GM13}). In the case of general potentials,
existence, uniqueness and  regularity of the weak solution have been shown in 
\cite{CF2} (see also \cite{FY} for an optimal control problem) by using an abstract approach. 
In the problem considered by Gal~\cite{G1} the Laplace--Beltrami term 
was missing in the third condition in \eqref{pier4} (thus, the boundary condition was of 
Wentzell type);  on the other hand, the presence  of the term  $- \DeltaG\muG $ 
actually enhances the dissipation mechanism in  \eqref{pier2}--\eqref{pier6}
and is helpful in order to recover a better regularity on the solution.
However, it is worth to point out that in \cite{CGS13} the presence of the convective term
$\nabla\rho\cdot u$ gives rise to further complications in the analysis.}

\pier{Some class of Cahn--Hilliard systems, possibly including dynamic boundary conditions, 
have collected an increasing interest in recent years: we can quote 
\cite{CFP, Kub12, MZ, PRZ, RZ, WZ} among other contributions. In case 
of no convective term in \eqref{pier2}, and assuming the homogeneous 
boundary condition $\dn \mu =0$ and the  condition
\eqref{pier5} with $\tauG>0$ and $\muG$ as a given datum, the problem 
has been first addressed in \cite{GiMiSchi}: the well-posedness and 
the large time behavior of solutions have been studied for regular 
potentials $f$ and $f_\Gamma$, as well as for various singular 
potentials like the ones in \eqref{logpot} and \eqref{obspot}. 
One can see \cite{GiMiSchi, GiMiSchi2}: in these two papers the authors 
were able to overcome the difficulties due to singularities using a set 
of assumptions for $f$ and $\fG $ that gives the 
role of the dominating potential to $f$ and entails some technical difficulties.  
The subsequent papers 
\cite{CGS3, CGS5, CGS4} follow a different approach, which
was firstly considered in \cite{CaCo} and \cite{CS} to investigate the Allen--Cahn 
equation with dynamic boundary conditions. This approach consists in 
letting $f_\Gamma$ be the leading potential with respect to $f$: 
by this, the analysis turns out to be simpler. In particular,
\cite{CGS3}~contains many results about existence, uniqueness 
and regularity of solutions for general potentials that include \accorpa{regpot}{logpot},
and are valid for both the viscous and pure cases, i.e., by assuming just $\tauO\geq0$.
The paper \cite{CF1} deals with the well-posedness of the same system, but in 
which also an additional mass constraint on the boundary is imposed. The recent contribution~\cite{LW}
deals with the physical derivation of some Cahn--Hilliard systems in the bulk and on the boundary,
arriving at the study of a model in which the two chemical potentials are completely independent. 
Finally, let us point out that the optimal control problems for \eqref{pier2}--\eqref{pier6} 
with the velocity as the control is thoroughly discussed in \cite{CGS14,CGS15}.}

The \pier{present} paper is organized as follows. 
In the next section, we list our assumptions and notations,
recall the properties already known 
and state our result on the \longtime\ \bhv.
The last section is devoted to the corresponding proof.


\section{Statement of the problem and results}
\label{STATEMENT}
\setcounter{equation}{0}

In this section, we state precise assumptions and notations and present our results.
First of all, the set $\Omega\subset\erre^3$ is assumed to be bounded, connected and smooth.
As in the introduction, $\nu$~is the outward unit normal vector field on $\Gamma:=\partial\Omega$, 
and $\dn$ and $\DeltaG$ stand for the corresponding normal derivative
and the Laplace--Beltrami operator, respectively.
Furthermore, we denote by $\nablaG$ the surface gradient
and write $|\Omega|$ and $|\Gamma|$ 
for the volume of $\Omega$ and the area of~$\Gamma$, respectively.
Moreover, we widely use the notations
\Beq
  Q_t := \Omega \times (0,T)
  \aand
  \Sigma_t := \Gamma \times (0,T)
  \quad \hbox{for $0<t\leq+\infty$}.
  \label{defQS}
\Eeq
Next, if $X$ is a Banach space, $\norma\cpto_X$ denotes both its norm and the norm of~$X^3$, 
and the symbols $X^*$ and $\<\cpto,\cpto>_X$ stand for the dual space of $X$
and the duality pairing between $X^*$ and~$X$.
The only exception from the convention for the norms is given
by the Lebesgue spaces~$L^p$, for $1\leq p\leq\infty$, 
whose norms will be denoted by~$\norma\cpto_p$. 
Furthermore, we~put
\Bsist
  && H := \Ldue \,, \quad  
  V := \Huno 
  \aand
  W := \Hdue,
  \label{defspaziO}
  \\
  && \HG := \LdueG \,, \quad 
  \VG := \HunoG 
  \aand
  \WG := \HdueG,
  \label{defspaziG}
  \\
  && \calH := H \times \HG \,, \quad
  \calV := \graffe{(v,\vG) \in V \times \VG : \ \vG = v\suG}
  \non
  \\
  && \aand
  \calW := \bigl( W \times \WG \bigr) \cap \calV \,.
  \label{defspaziprod}
\Esist
In the following, we work in the framework of the Hilbert triplet
$(\calV,\calH,\calVp)$.
Thus, we have
$$\<(g,\gG),(v,\vG)>_{\calV}=\iO gv+\iG\gG\vG
\mbox{ \,for every $(g,\gG)\in\calH$ and $(v,\vG)\in\calV$.} $$

Now, we list our assumptions.
For the structure of our system, we postulate the following properties, 
which are slightly stronger than \juerg{those requested in} \cite{CGS13}:
\begin{align}
  & \hbox{$\tauO$ and $\tauG$ \ are strictly positive real numbers.}
  \label{hptau}
  \\
  & \Beta,\, \BetaG : \erre \to [0,+\infty]
  \quad \hbox{are convex, proper, and l.s.c., with} \quad
  \Beta(0) = \BetaG(0) = 0.
  \qquad
  \label{hpBeta}
  \\
  \separa
  & \Pi,\, \PiG : \erre \to \erre
  \quad \hbox{are of class $C^2$ with \Lip\ continuous first derivatives.}
  \qquad
  \label{hpPi}
  \\ 
  & \hbox{The functions $f:=\Beta+\Pi$ and $\fG:=\BetaG+\PiG$ are bounded from below}.
  \label{hpf}
\end{align}
Thus, in contrast to~\cite{CGS13},
the constants $\tauO$ and $\tauG$ are strictly positive, here,
and \eqref{hpf} holds in addition.
However, we remark that \juerg{this} assumption is \pier{very reasonable and} 
fulfilled by all of the potentials \accorpa{regpot}{obspot}.
We set, for convenience,
\Beq
  \beta := \partial\Beta \,, \quad
  \betaG := \partial\BetaG \,, \quad
  \pi := \Pi'\,,
  \aand
  \piG := \PiG',
  \label{defbetapi}  
\Eeq
and assume that, with some positive constants $C$ and $\eta$,
\Beq
  D(\betaG) \subseteq D(\beta)
  \aand
  |\betaz(r)| \leq \eta |\betaGz(r)| + C
  \quad \hbox{for every $r\in D(\betaG)$}.
  \label{hpCC}
\Eeq
\Accorpa\HPstruttura hptau hpCC
In \eqref{hpCC}, the symbols $D(\beta)$ and $D(\betaG)$ 
denote the domains of $\beta$ and~$\betaG$, respectively.
More generally, we use the notation $D(\calG)$ 
for every maximal monotone graph $\calG$ in $\erre\times\erre$,
as well as for the \juerg{associated} maximal monotone operators induced on $L^2$ spaces.
Moreover, for $r\in D(\calG)$,
$\calG^\circ(r)$ stands for the element of $\calG(r)$ having minimum modulus.
\juerg{
\Brem
Notice that, physically speaking, the compatibility condition \eqref{hpCC} means that
the thermodynamic force driving the phase separation on the surface is stronger than
the one in the bulk.  
\Erem
}
For the data, we make the following assumptions:
\Bsist
  && u \in (\LL2{\Lx3})^3 
  \quad \hbox{with} \quad 
  \dt u \in (\LL2{\Lx{3/2}})^3\,.
  \label{hpu}
  \\
  && \div u = 0 \quad \hbox{in $\Qi$}
  \aand
  u \cdot \nu = 0 \quad \hbox{on $\Si$}\,.
  \label{hpubis}
  \\
  && (\rhoz \,,\, \rhoGz) \in \calW \,, \quad
  \betaz(\rhoz) \in H
  \aand
  \betaGz(\rhoGz) \in \HG \,.
  \qquad
  \label{hprhoz}
  \\
  && \mz := \frac { \iO\rhoz + \iG\rhoGz } { |\Omega| + |\Gamma| }
  \quad \hbox{belongs to the interior of $D(\betaG)$}\, .
  \label{hpmz}
\Esist
\Accorpa\HPdati hpu hpmz
Since we are assuming~\eqref{hptau} and \HPdati,
the regularity level \juerg{required for} the notion of solution on a finite time interval
is higher than \juerg{the one} in~\cite{CGS13}.
Namely, a~solution on $(0,T)$ is a triple of pairs $\jsoluz$ that satisfies
\Bsist
  && \Mu \in \L\infty\calW, 
  \label{regmu}
  \\
  && \Rho \in \W{1,\infty}\calH \cap \H1\calV \cap \L\infty\calW, 
  \label{regrho}
  \\
  && \Zeta \in \L\infty\calH .
  \label{regzeta}
\Esist
\Accorpa\Regsoluz regmu regzeta
However, we write $\soluz$ instead of $\jsoluz$, in order to simplify the notation.
As far as the problem under study is concerned, 
we still state it in a weak form as in~\cite{CGS13}, 
\pier{on account of} the assumptions \eqref{hpubis} on~$u$.
Namely, we require that
\Bsist
  && \iO \dt\rho \, v
  + \iG \dt\rhoG \, \vG
  - \iO \rho u \cdot \nabla v
  + \iO \nabla\mu \cdot \nabla v
  + \iG \nablaG\muG \cdot \nablaG\vG
  = 0
  \non
  \\[1mm]
  && \quad \hbox{\aet\ and for every $(v,\vG)\in\calV$},
  \label{prima}
  \\[2mm]
  \separa
  && \tauO \iO \dt\rho \, v
  + \tauG \iG \dt\rhoG \, \vG
  + \iO \nabla\rho \cdot \nabla v
  + \iG \nablaG\rhoG \cdot \nablaG\vG
  \non
  \\[1mm]
  && \quad {}
  + \iO \bigl( \zeta + \pi(\rho) \bigr) v
  + \iG \bigl( \zetaG + \piG(\rhoG) \bigr) \vG
  = \iO \mu v 
  + \iG \muG \vG
  \non
  \\[1mm]
  && \quad \hbox{\aet\ and for every $(v,\vG)\in\calV$},
  \label{seconda}
  \\[2mm]
  && \zeta \in \beta(\rho) \quad \aeQ
  \aand
  \zetaG \in \betaG(\rhoG) \quad \aeS,
  \label{terza}
  \\[2mm]
  && \rho(0) = \rhoz
  \quad \aeO \,.
  \label{cauchy}
\Esist
\Accorpa\Pbl prima cauchy
However, every solution also satisfies the boundary value problem presented in the introduction.
The basic well-posedness and regularity results are given by~\cite[Thms.~2.3 and~2.6]{CGS13}.
We collect them in the following

\Bthm
\label{RecallCGS}
Assume \HPstruttura\ for the structure and \HPdati\ for the data,
and let $T\in(0,+\infty)$.
Then  problem \Pbl\ has at least one solution $\soluz$ satisfying \Regsoluz.
Moreover, the solution is unique if 
at least one of the operators $\beta$ and $\betaG$ is single-valued.
\Ethm

We obviously deduce the following consequence:

\Bcor
\label{Globalsolution}
In addition to \HPstruttura,
assume that at least one of the operators $\beta$ and $\betaG$ is single-valued.
Moreover, assume \HPdati\ for the data.
Then there exists a unique 6-tuple $\soluz$ defined on $(0,+\infty)$
that fulfils \Regsoluz\ and solves \Pbl\ for every $T\in(0,+\infty)$.
\Ecor

At this point, given a solution $\soluz$, 
our aim is investigating its \longtime\ \bhv,
namely, the \omegalimit\
(which we simply term~$\omega$ for brevity)
of the component~$\Rho$.
We notice that \pier{the property} \eqref{regrho} \pier{holds} 
for every $T\in(0,+\infty)$\pier{, which}
implies that $\Rho$ belongs to~$\CC0\calV$\pier{.
Hence,} the next definition is meaningful.
We~set
\begin{align}
   \omega := \bigl\{
  \Rhoo=\displaystyle\limn\Rho(\tn)
  \quad &\hbox{in the weak topology of $\calV$}
  \non
  \\
  &\hbox{for some sequence $\{\tn\nearrow+\infty\}$}
  \bigr\}\,,
  \label{omegalim}
\end{align}
and look for the relationship between $\omega$ and the set of stationary solutions
to the system \juerg{obtained} from \accorpa{prima}{terza} by ignoring the convective term.
Indeed, assumption \eqref{hpu} implies that 
\gianni{%
\Beq
  u \in \pier{(}\HH1{\Lx{3/2}}\pier{)^3\,,}
  \label{utozero}
\Eeq
whence $u(t)$ tends to zero strongly in $\Lx{3/2}$ as $t$ tends to infinity.}
It is immediately seen from \eqref{prima} 
that the components $\mu$ and $\muG$
of every stationary solution are constant functions
and that the constant values they assume are the same.
Therefore, by a stationary solution we mean 
a quadruplet $\soluzs$ satisfying for some $\mus\in\erre$ the conditions
\Bsist
  && \Rhos \in \calV
  \aand
  \Zetas \in \calH,
  \label{regs}
  \\[2mm]
  && \iO \nabla\rhos \cdot \nabla v
  + \iG \nablaG\rhoGs \cdot \nablaG\vG
  + \iO \bigl( \zetas + \pi(\rhos) \bigr) v
  + \iG \bigl( \zetaGs + \piG(\rhoGs) \bigr) \vG
  \non
  \\[1mm]
  && \quad {}
  = \iO \mus v 
  + \iG \mus \vG
  \quad \hbox{for every $(v,\vG)\in\calV$},
    \label{secondas}
  \\[2mm]
  && \zetas \in \beta(\rhos) \quad \aeO
  \aand
  \zetaGs \in \betaG(\rhoGs) \quad \aeG \,.
  \label{terzas}
\Esist
It is not difficult to show that \accorpa{regs}{secondas}
imply that the pair $\Rhos$ belongs to $\calW$ and 
satisfies the boundary value problem
\Bsist
  && - \Delta\rhos + \zetas + \pi(\rhos) = \mus
  \quad \aeO\,,
  \non
  \\[1mm]
  && \dn\rhos - \DeltaG\rhoGs + \zetaGs + \piG(\rhoGs) = \mus
  \quad \aeG \,.
  \non
\Esist
Here is our result:

\Bthm
\label{Omegalimit}
Let the assumptions of Corollary~\ref{Globalsolution} be satisfied,
and let $\soluz$ be the unique global solution on $(0,+\infty)$.
Then the \omegalimit\ \eqref{omegalim} is nonempty.
Moreover, for every $\Rhoo\in\omega$,
there exist $\mus\in\erre$ and a solution $(\rhos,\rhoGs,\zetas,\zetaGs)$ to \accorpa{regs}{terzas}
such that $\Rhoo=\Rhos$.
\Ethm


\section{Auxiliary material}
\label{AUX}
\setcounter{equation}{0}

The proof of Theorem~\ref{Omegalimit} will be performed in the last section.
Our argument needs some tools which we collect in the present section.
In particular, we use the generalized mean value, the related spaces 
and the operator~$\calN$ which we introduce now.
However, we proceed very shortly and refer to \cite[Sect.~2]{CGS13}
for further details.
We set
\Beq
  \mean\gstar := \frac { \< \gstar,(1,1) >_{\calV} } { |\Omega| + |\Gamma| }
  \quad \hbox{for $\gstar\in\calVp$}
  \label{genmean}
\Eeq
and observe~that
\Beq
  \mean(v,\vG) = \frac { \iO v + \iG \vG } { |\Omega| + |\Gamma| }
  \quad \hbox{if $\gstar=(v,\vG)\in\calH$} \,.
  \label{usualmean}
\Eeq
Notice that the constant $\mz$ appearing in our assumption \eqref{hpmz} 
is nothing but  the mean value $\mean(\rhoz,\rhoGz)$,
and that taking $(v,\vG)=(|\Omega|+|\Gamma|)^{-1}(1,1)$ in \pier{\eqref{prima}}
yields the conservation property for the component $\Rho$ of the solution,
\Beq
  \dt \mean\Rho = 0,
  \quad \hbox{whence} \quad
  \mean \Rho(t) = \mz
  \quad \hbox{for every $t\in[0,T]$}.
  \label{conservation}
\Eeq
We also stress that the function
\Beq
  \calV \ni (v,\vG) \mapsto
  \bigl( \norma{\nabla v}_{\L2H}^2 + \norma{\nablaG\vG}_{\L2\HG}^2 + |\mean(v,\vG)|^2 \bigr)^{1/2}
  \label{normacalV}
\Eeq
yields a Hilbert norm on $\calV$ which is equivalent to the natural one.
Now, we set
\Beq
  \calVsz := \graffe{ \gstar\in\calVp : \ \mean\gstar = 0 } , \quad
  \calHz := \calH \cap \calVsz,
  \aand
  \calVz := \calV \cap \calVsz, 
  \label{defcalVz}
\Eeq
\Accorpa\Defspazi defspaziO defcalVz
and notice that the function
\Beq
  \calVz \ni (v,\vG) \mapsto \norma{(v,\vG)}_{\calVz}
  := \bigl( \norma{\nabla v}_{\L2H}^2 + \norma{\nablaG\vG}_{\L2\HG}^2 \bigr)^{1/2}
  \label{normaVz}
\Eeq
is a Hilbert norm on $\calVz$ which is equivalent to the usual one.
Next, we define the operator $\calN:\calVsz\to\calVz$ 
(which~will be applied to $\calVsz$-valued functions as well)
as~follows.
For every element $\gstar\in\calVsz$,
\Bsist
  && \hbox{$\calN\gstar=(\calNO\gstar,\calNG\gstar)$ is the unique pair $\Xi\in\calVz$ such that}
  \non
  \\
  \noalign{\smallskip}
  && \iO \nabla\xi \cdot \nabla v + \iG \nablaG\xiG \cdot \nablaG\vG
  = \< \gstar , (v,\vG) >_{\calV}
  \quad \hbox{for every $(v,\vG)\in\calV$}.
  \qquad
  \label{defN}
\Esist
It turns out that $\calN$ is well defined, linear, symmetric, and bijective.
Therefore, if we set
\Beq
  \norma\gstar_* := \norma{\calN\gstar}_{\calVz},
  \quad \hbox{for $\gstar\in\calVsz$},
  \label{normastar}
\Eeq
then we obtain a Hilbert norm on $\calVsz$
(equivalent to the norm induced by the norm of~$\calVp$).
Furthermore, we notice that
\Beq
  \< \dt\gstar , \calN\gstar >_{\calV}
  =  \frac 12 \, \frac d{dt} \, \norma\gstar_*^2
  \quad \aet,
  \quad \hbox{ for every $\gstar \in \H1\calVsz$}.
  \label{propNdta}
\Eeq
Finally, it is easy to see that $\calN\gstar$ belongs to $\calW$ whenever $\gstar\in\calHz$, 
and that
\Beq
  \normaWW{\calN\gstar} \leq \CO \normaHH\gstar
  \quad \hbox{for every $\gstar\in\calHz$},
  \label{regN}
\Eeq
where $\CO$ depends only on $\Omega$.

In performing our estimates, we will repeatedly use the Young inequality
\Beq
  a\,b \leq \delta\,a^2 + \frac 1{4\delta} \, b^2
  \quad \hbox{for all $a,b\in\erre$ and $\delta>0$},
  \label{young}
\Eeq
as well as H\"older's inequality and the Sobolev inequality
\Beq
  \norma v_p
  \leq \CO \normaV v
  \quad \hbox{for every $p\in[1,6]$ and $v\in V$},
  \label{sobolev}
\Eeq
\juerg{which is} related to the continuous embedding $V\subset L^p(\Omega)$
for $p\in[1,6]$
(since $\Omega$ is three-dimensional, bounded and smooth).
In~particular, by also using the equivalent norm \eqref{normacalV} on~$\calV$,
we have that
\Beq
  \norma v_6^2
  \leq \CO \bigl( \norma{\nabla v}_{\L2H}^2 + \norma{\nablaG\vG}_{\L2\HG}^2 + |\mean(v,\vG)|^2 \bigr)
  \label{sobolevbis}
\Eeq
for every $(v,\vG)\in\calV$.
In both \eqref{sobolev} and \eqref{sobolevbis}, the constant $\CO$ depends only on~$\Omega$.
Furthermore, we will owe to a well-known fact from interpolation theory.
\gianni{By \cite[Thm.~5.2.1\pier{,} p.~109]{BL}, we have that 
\Beq
  \bigl( \Lx3 , \Lx{3/2} \bigr)_{1/2,2}
  = \Lx2
  \non
\Eeq
where the interpolation can be understood also in the sense of the trace method
(see \cite[Sect.~3.12]{BL} for the equivalence between various interpolation methods).
It follows that there hold} the continuous embedding and the related inequality
\Bsist
  && \HH1{\Lx{3/2}} \cap \LL2{\Lx3} \subset \LL\infty\Ldue,
  \non
  \\[3pt]
  && \norma u_{\LL\infty\Ldue}
  \leq \CO \, \bigl(
    \norma u_{\LL2{\Lx3}}
    + \norma{\dt u}_{\LL2{\Lx{3/2}}}
  \bigr)
  \non
  \\
  && \quad \hbox{for every $v\in\HH1{\Lx{3/2}}\cap\LL2{\Lx3}$},
  \label{interpolaz}
\Esist
where $\CO$ depends only on~$\Omega$.

Finally, as far as constants are concerned,
we \juerg{employ the following general rule}: the small-case symbol $c$ stands 
for different constants which depend only
on~$\Omega$, the structure of our system 
and the norms of the data involved in the assumptions \HPdati.
A~notation like~$c_\delta$ (in particular, with $\delta=T$) 
allows the constant to depend on the positive
parameter~$\delta$, in addition.
Hence, the meaning of $c$ and $c_\delta$ may
change from line to line and even within the same chain of inequalities.  
On the contrary, we mark the constants that we want to refer~to
by using a different notation (e.g., a~capital letter).


\section{\Longtime\ \bhv}
\label{PROOFS}
\setcounter{equation}{0}

This section is devoted to the proof of Theorem~\ref{Omegalimit}.
Thus, we fix any global solution $\soluz$ once and for all.
Our arguments relies on some global a~priori estimates
and on the study of the \bhv\ of such a solution on intervals 
of a fixed length~$T$ whose \juerg{endpoints} tend to infinity.

To keep the \juerg{paper at a reasonable length},
we often proceed formally.
However, we notice that the estimates \juerg{to be  obtained in this way
can be performed rigorously by}
acting on the solution to a proper regularizing or discrete problem.
Indeed, the solution found in \cite{CGS13}
\juerg{was} constructed in this way:
first, a~regularized $\eps$-problem was introduced 
by replacing the graphs $\beta$ and $\betaG$ 
by their Yosida regularizations $\betaeps$ and $\betaGeps$.	
However, the same argument \juerg{would work if $\betaeps$ and $\betaGeps$
were} smooth approximation of $\beta$ and $\betaG$
with analogous boundedness and convergence properties
(like the $C^\infty$ approximations introduced in \cite[Sect.~3]{GiRo}).
In order to solve such an approximating problem,
a~Faedo-Galerkin scheme depending on a parameter $n\in\enne$ \juerg{can be} used.
Its solution is smooth according to the smoothness of the nonlinearities that appear in the $\eps$-problem.
Then, the solution to the original problem is constructed
by first letting $n$ tend to infinity
and then letting $\eps$ tend to zero.
Thus, our procedure would be completely rigorous 
if \juerg{it were} performed on one of the above approximating solutions.
Indeed, \juerg{the estimates established in this way would be} uniform with respect to the parameters involved.
We now start proving such global estimates.

\step
First global estimate

We write the equations \eqref{prima} and \eqref{seconda} at the time~$s$
and test them by $\Mu(s)\in\calV$ and $\dt\Rho(s)\in\calV$, respectively.
Then, we integrate with respect to $s$ over $(0,t)$ with an arbitrary~$t>0$, 
sum up and rearrange.
We~obtain
\Bsist
  && \intQt \dt\rho \, \mu
  + \intSt \dt\rhoG \, \muG
  + \intQt |\nabla\mu|^2
  + \intSt |\nablaG\muG|^2  
  + \tauO \intQt |\dt\rho|^2
  + \tauG \intSt |\dt\rhoG|^2
  \non
  \\
  && \quad {}
  + \frac 12 \iO |\nabla\rho(t)|^2
  + \frac 12 \iG |\nablaG\rhoG(t)|^2
  + \iO f(\rho(t))
  + \iG \fG(\rhoG(t))
  \non
  \\
  && = \intQt \rho u \cdot \nabla\mu
  + \intQt \mu \dt\rho
  + \intSt \muG \dt\rhoG
  \non
  \\
  && \quad {}
  + \frac 12 \iO |\nabla\rhoz|^2
  + \frac 12 \iG |\nablaG\rhoGz|^2
  + \iO f(\rhoz)
  + \iG \fG(\rhoGz).
  \non
\Esist
Four integrals obviously cancel out,
the ones on the \lhs\ containing $f$ and $\fG$ are bounded from below by~\eqref{hpf},
and the terms on the \rhs\ involving the initial values are finite by~\eqref{hprhoz}.
We deal with the convective term \pier{using} the Young and \Holder\ inequalities,
the Sobolev type inequality \eqref{sobolevbis}, and the conservation property~\eqref{conservation}.
We~have
\Bsist
  && \intQt \rho u \cdot \nabla\mu
  \leq \frac 12 \intQt |\nabla\mu|^2
  + \frac 12 \iot \norma{u(s)}_3^2 \, \norma{\rho(s)}_6^2 \, ds
  \non 
  \\
  && \leq \frac 12 \intQt |\nabla\mu|^2
  + c \iot \norma{u(s)}_3^2 \, \bigl( \juerg{\|\nabla\rho(s)\|_{H}^2 + \|\nablaG\rhoG(s)\|_{\HG}^2} + \mz^2 \bigr) \, ds .
  \non
\Esist
Since the function $s\mapsto\norma{u(s)}_3^2$ belongs to $L^1(0,+\infty)$ by~\eqref{hpu},
we can apply the Gronwall lemma on $(0,+\infty)$ and obtain~that
\Bsist
  && \intQi |\nabla\mu|^2
  + \intSi |\nablaG\muG|^2
  + \intQi |\dt\rho|^2
  + \intSi |\dt\rhoG|^2
  < + \infty\,,
  \label{globint}
  \\[3pt]
  && f(\rho) \in \LL\infty\Luno
  \aand
  \fG(\rhoG) \in \LL\infty\LunoG\,,
  \label{stimef}
\Esist
\juerg{as well as} $(\nabla\rho,\nablaG\rhoG)\in(\LL\infty\calH)^3$.
From this, by accounting for the conservation property~\eqref{conservation} once more,
we conclude that
\Beq
  \Rho \in \LL\infty\calV .
  \label{primastima}
\Eeq

\step
Consequence

By using the quadratic growth of $\Pi$ and $\PiG$ 
implied by the \Lip\ continuity of their derivatives,
and combining with \eqref{stimef} with~\eqref{primastima},
we deduce that
\Beq
  \Beta(\rho) \in  \LL\infty\Luno
  \aand
  \BetaG(\rhoG) \in \LL\infty\LunoG \,.
  \label{stimeBeta}
\Eeq

\step
Second global estimate

We formally differentiate the equations \eqref{prima} and \eqref{seconda}
with respect to time\juerg{, where we argue as if $\beta$ and $\betaG$ were smooth functions
and write} $\beta(\rho)$ and $\betaG(\rhoG)$
instead of $\zeta$ and~$\zetaG$ (see~\eqref{terza}).
We obtain that
\Bsist
  && \iO \dt^2\rho \, v
  + \iG \dt^2\rhoG \, \vG
  + \iO \nabla\dt\mu \cdot \nabla v
  + \iG \nablaG\dt\muG \cdot \nablaG\vG
  \non
  \\
  && = \iO \bigl( \dt\rho \, u + \rho \dt u \bigr) \cdot \nabla v
  \label{dtprima}
  \\
  \separa
  && \tauO \iO \dt^2\rho \, v
  + \tauG \iG \dt^2\rhoG \, \vG
  + \iO \nabla\dt\rho \cdot \nabla v
  + \iG \nablaG\dt\rhoG \cdot \nablaG\vG
  \non
  \\
  && \quad {}
  + \iO \beta'(\rho) \dt\rho \, v
  + \iG \betaG'(\rhoG) \dt\rhoG \, \vG
  \non
  \\
  && = \iO \dt\mu \, v 
  + \iG \dt\muG \, \vG
  - \iO \pi'(\rho) \dt\rho \, v
  - \iG \piG'(\rhoG) \dt\rhoG \, \vG 
  \label{dtseconda}
\Esist
\Aet\ and for every $(v,\vG)\in\calV$.
By recalling that $\dt\Rho$ is $\calVz$-valued by \eqref{conservation},
so that $\calN\dt\Rho$ is well defined,
we write the above equations at the time~$s$
and test them by $\calN\dt\Rho(s)$ and $\dt\Rho(s)$, respectively.
Then, we integrate with respect to $s$ over $(0,t)$ with an arbitrary $t>0$ and sum~up.
We~have
\Bsist
  && \iot \< \dt^2\Rho(s) , \calN\dt\Rho(s) >_\calV \, ds
  \non
  \\
  && \quad {}
  + \intQt \nabla\dt\mu \cdot \nabla\calNO(\dt\Rho)
  + \intSt \nablaG\dt\muG \cdot \nablaG\calNG(\dt\Rho)
  \non
  \\
  && \quad {}
  + \frac \tauO 2 \iO |\dt\rho(t)|^2
  + \frac \tauG 2 \iG |\dt\rhoG(t)|^2
  + \intQt |\nabla\dt\rho|^2
  + \intSt |\nablaG\dt\rhoG|^2
  \non
  \\
  && \quad {}
  + \intQt \beta'(\rho) |\dt\rho|^2
  + \intSt \betaG'(\rhoG) |\dt\rhoG|^2
  \non
  \\
  \separa
  && = \intQt \bigl( \dt\rho \, u + \rho \dt u \bigr) \cdot \nabla\calNO(\dt\Rho)\,
  \juerg{+ \,\frac \tauO 2 \iO |\dt\rho(0)|^2 + \frac \tauG 2 \iG |\dt\rhoG(0)|^2} \non
  \\
  && \quad {}
  + \intQt \dt\mu \, \dt\rho 
  + \intSt \dt\muG \, \dt\juerg{\rhoG}
  - \intQt \pi'(\rho) |\dt\rho|^2
  - \intSt \piG'(\rhoG) |\dt\rhoG|^2 .
  \non
\Esist
The integrals containing $\dt\mu$ and $\dt\muG$ cancel out
by the definition \eqref{defN} of~$\calN$
(with the choices $\gstar=\dt\Rho(s)$ and $(v,\vG)=\dt\Mu(s)$),
and the terms involving $\beta'$ and $\betaG'$ are nonnegative.
Moreover, the last two integrals on the \rhs\ are bounded
by \eqref{globint} and the \Lip\ continuity of $\pi$ and~$\piG$.
Therefore, \pier{owing to} \eqref{propNdta} for the first term on the \lhs,
we deduce that
\Bsist
  && \frac 12 \, \norma{\dt\Rho(t)}_*^2
  + \frac \tauO 2 \iO |\dt\rho(t)|^2
  + \frac \tauG 2 \iG |\dt\rhoG(t)|^2
  + \intQt |\nabla\dt\rho|^2
  + \intSt |\nablaG\dt\rhoG|^2
  \non
  \\
  && \leq \frac 12 \, \norma{\dt\Rho(0)}_*^2
  + \frac \tauO 2 \iO |\dt\rho(0)|^2
  + \frac \tauG 2 \iG |\dt\rhoG(0)|^2
  \non
  \\
  && \quad {}
  + \intQt \dt\rho \, u \cdot \nabla\calNO(\dt\Rho)
  + \intQt \rho \dt u \cdot \nabla\calNO(\dt\Rho)
  + c \,.
  \label{persecondastima}
\Esist
Thus, is suffices to obtain a bound for the $\calH$-norm of $\dt\Rho(0)$ and to estimate the last two integrals.
For the first aim, we write the equations \eqref{prima} and \eqref{seconda} at the time $t=0$,
test them by $\Mu(0)$ and $\dt\Rho(0)$, respectively, and sum~up.
Then, we account for the regularity of $\rhoz$ ensured by \eqref{hprhoz}
and integrate by parts the terms involving $\nabla\rhoz$ and~$\nablaG\rhoGz$.
We obtain
\Bsist
  && \iO \dt\rho(0) \, \mu(0) 
  + \iG \dt\rhoG(0) \, \muG(0)
  + \iO |\nabla\mu(0)|^2
  + \iG |\nablaG\muG(0)|^2
  \non
  \\
  && \quad {}
  + \tauO \iO |\dt\rho(0)|^2
  + \tauG \iG |\dt\rhoG(0)|^2
  - \iO \Delta\rhoz \, \dt\rho(0)
  - \iG \DeltaG\rhoGz \, \dt\rhoG(0)
  \non
  \\
  && \quad {}
  + \iO (\beta + \pi)(\rhoz) \dt\rho(0)
  + \iG (\betaG + \piG)(\rhoGz) \dt\rhoG(0)
  \non
  \\
  \separa
  && = \iO \rhoz u(0) \cdot \nabla\mu(0)
  + \iO \mu(0) \dt\rho(0)
  + \iG \muG(0) \dt\rhoG(0).
  \non
\Esist
Four integrals obviously cancel each other.
Now, we rearrange and use the Young inequality and the \Lip\ continuity of $\pi$ and~$\piG$.
\pier{In view of} the full~\eqref{hprhoz}, we infer~that
\Bsist
  && \iO |\nabla\mu(0)|^2
  + \iG |\nablaG\muG(0)|^2 
  + \frac \tauO 2 \iO |\dt\rho(0)|^2
  + \frac \tauG 2 \iG |\dt\rhoG(0)|^2 
  \non
  \\
  && \leq c \bigl( \normaWW\Rhoz^2 + 1 + \normaH{\beta(\rhoz)}^2 + \normaHG{\betaG(\rhoGz)}^2 \bigr)
  + \iO \rhoz u(0) \cdot \nabla\mu(0)
  \non
  \\
  && \leq c + \iO \rhoz u(0) \cdot \nabla\mu(0).
  \non
\Esist
In order to estimate the last integral,
we use the \Holder\ and Young inequalities,
the continuous embedding $W\subset\Lx\infty$, and the interpolation inequality \eqref{interpolaz},
\juerg{to conclude that}
\Bsist
  && \iO \rhoz u(0) \cdot \nabla\mu(0)
  \leq \norma\rhoz_\infty \, \pier{\norma {u(0)}_2} \, \norma{\nabla\mu(0)}_2
  \non
  \\
  && \leq \frac 12 \iO |\nabla\mu(0)|^2
  + c \, \normaW\rhoz^2 \bigl( \pier{\norma u_{\LL2{\Lx3}}^2 + \norma{\dt u}_{\LL2{\Lx{3/2}}}^2} \bigr)\,,  
  \non
\Esist
and we observe that all of the last norms are finite \juerg{by virtue of the} 
assumptions \eqref{hprhoz} and \eqref{hpu} on $\rhoz$ and~$u$.
By combining this and the above inequality, 
we obtain the desired bound for the $\calH$-norm of~$\dt\Rho(0)$.

Finally, we estimate the last two integrals on the \rhs\ of~\eqref{persecondastima}
by using the \Holder, Sobolev and Young inequalities (in particular~\eqref{sobolevbis}),
the interpolation inequality \eqref{interpolaz},
the conservation property~\eqref{conservation},
the regularity inequality~\eqref{regN} for~$\calN$,
and the already established estimate \eqref{primastima}.
We~have that
\Bsist
  && \intQt \dt\rho \, u \cdot \nabla\calNO(\dt\Rho)
  \leq \iot \norma{\dt\rho(s)}_4 \, \norma{u(s)}_2 \, \norma{\nabla\calNO(\dt\Rho(s))}_4 \, ds
  \non
  \\
  && \leq c \, \norma u_{\LL\infty H} \iot \normaV{\dt\rho(s)} \, \normaWW{\calN(\dt\Rho(s))} \, ds
  \non
  \\
  && \leq \frac 12 \, \Bigl( \intQt |\nabla\dt\rho|^2 + \intSt |\nablaG\dt\rhoG|^2 \Bigr)
  + c \iot \normaHH{\dt\Rho(s)}^2 \, ds
  \non
  \\
  && \leq \frac 12 \, \Bigl( \intQt |\nabla\dt\rho|^2 + \intSt |\nablaG\dt\rhoG|^2 \Bigr)
  + c\,,
  \non
\Esist
as well as, with a similar argument,
\Bsist
  && \intQt \rho \dt u \cdot \nabla\calNO(\dt\Rho)
  \leq \iot \norma{\rho(s)}_6 \, \norma{\dt u(s)}_{3/2} \, \norma{\nabla\calNO(\dt\Rho(s))}_6 \, ds
  \non
  \\
  && \leq c \, \norma\rho_{\LL\infty V} \, \norma{\dt u}_{\LL2{\Lx{3/2}}} \, \norma{\dt\Rho}_{\LL2\calH}
  \leq c \,.
  \non
\Esist
Coming back to \eqref{persecondastima} 
and taking these estimates into account, we conclude that
\Beq
  \dt\Rho \in \LL\infty\calH\cap\LL2\calV \,.
  \label{secondastima}
\Eeq

\step
Third global estimate

Also in this step, we \juerg{argue as if} the graphs $\beta$ and $\betaG$
were smooth functions
(in particular, we write $\zeta=\beta(\rho)$ and $\zetaG=\beta(\rhoG)$)
and first notice that the inclusion $D(\betaG)\subseteq D(\beta)$ (see~\eqref{hpCC}) 
and assumption \eqref{hpmz} imply that
\Beq
  \beta(r) (r-\mz)
  \geq \delta_0 |\beta(r)| - C_0
  \aand
  \betaG(r) (r-\mz)
  \geq \delta_0 |\betaG(r)| - C_0
  \label{trickMZ}
\Eeq
for every $r$ belonging to the respective domains,
where $\delta_0$ and $C_0$ are some positive constants that depend only on $\beta$, $\betaG$
and on the position of $\mz$ in the interior of~$D(\betaG)$ and of~$D(\beta)$
(see, e.g., \cite[p.~908]{GiMiSchi}).
Now, we recall the conservation property \eqref{conservation}
and test \eqref{prima} and \eqref{seconda} by 
$\calN(\rho-\mz,\rhoG-\mz)$ and $(\rho-\mz,\rhoG-\mz)$, respectively.
Then, we sum up without integrating with respect to time.
We obtain \juerg{that, almost everywhere in $(0,+\infty)$,}
\Bsist
  && \< \dt\Rho , \calN(\rho-\mz,\rhoG-\mz) >_{\calV} 
  \non
  \\[1mm]
  && \quad {}
  + \iO \nabla\mu \cdot \nabla\calNO(\rho-\mz,\rhoG-\mz)
  + \iG \nablaG\muG \cdot \nablaG\calNG(\rho-\mz,\rhoG-\mz)
  \non
  \\
  && \quad {}
  + \tauO \iO \dt\rho \, (\rho-\mz)
  + \tauG \iG \dt\rhoG \, (\rhoG-\mz)
  + \iO |\nabla\rho|^2
  + \iG |\nablaG\rhoG|^2
  \non
  \\
  && \quad {}
  + \iO \beta(\rho) (\rho-\mz)
  + \iG \betaG(\rhoG) (\rhoG-\mz)
  \non
  \\
  && = \iO \rho u \cdot \nabla(\calNO(\rho-\rhoz,\rhoG-\rhoGz))
  + \iO \mu (\rho-\mz) 
  + \iG \muG (\rhoG-\mz)
  \non
  \\
  && \quad {}
  - \iO \pi(\rho) (\rho-\mz)
  - \iG \piG(\rhoG) (\rhoG-\mz)\,.
  \non
\Esist
All of the integrals involving $\mu$ and $\muG$ cancel out by \eqref{defN}.
Now, we owe to~\eqref{trickMZ}, keep just the positive contributions on the \lhs\
and move the other terms \pier{on} the \rhs.
By also accounting for the \Lip\ continuity of $\pi$ and~$\piG$,
the Sobolev inequality related to the continuous embedding $V\subset\Lx4$, 
the interpolation inequality \eqref{interpolaz} on~$u$,
the regularity inequality \eqref{regN} for~$\calN$,
\eqref{primastima} and~\eqref{secondastima}, we deduce~that
\Bsist
  &&  \iO |\nabla\rho|^2
  + \iG |\nablaG\rhoG|^2
  + \delta_0 \iO |\beta(\rho)|
  + \delta_0 \iG |\betaG(\rhoG)|
  \non
  \\[2mm]
  && \leq \norma{\dt\Rho}_{\LL\infty\calVzp} \, \norma{\calN(\rho-\mz,\rhoG-\mz)}_{\LL\infty\calVz}
  \non
  \\[3pt]
  && \quad {}
  + c \, \norma{\dt\Rho}_{\LL\infty\calH} \, \norma{(\rho-\mz,\rhoG-\mz)}_{\LL\infty\calH}
  \non
  \\[3pt]
  && \quad {}
  + c \, \norma{(\pi(\rho),\piG(\rhoG))}_{\LL\infty\calH} \, \norma{(\rho-\mz,\rhoG-\mz)}_{\LL\infty\calH}
  \non
  \\
  && \quad {}
  + \norma\rho_{\LL\infty{\Lx4}} \, \norma u_{\LL\infty\Ldue} \, \norma{\nabla(\calNO(\rho-\rhoz,\rhoG-\rhoGz))}_{\LL\infty{\Lx4}}
  \non
  \\[1mm]
  && \leq c
  + c \, \norma{\calN(\rho-\rhoz,\rhoG-\rhoGz)}_{\LL\infty\calW}
  \leq c 
  + c \, \norma{(\rho-\rhoz,\rhoG-\rhoGz)}_{\LL\infty\calH}
  \leq c \,.
  \non
\Esist
Since this holds \Aet,
we have (in particular) that
\Beq
  \zeta \in \LL\infty\Luno
  \aand
  \zetaG \in \LL\infty\LunoG .
  \non
\Eeq
Now, we test \eqref{seconda} by $(1,1)$
and obtain \Aet
\Beq
  (|\Omega|+|\Gamma|) \mean\Mu 
  = \tauO \iO \dt\rho + \tauG \iG \dt\rhoG
  + \iO \zeta + \iG \zetaG
  + \iO \pi(\rho) + \iG \piG(\rhoG).
  \non
\Eeq
Thus, we infer that
\Beq
  \mean\Mu \in L^\infty(0,+\infty) \,.
  \label{terzastima}
\Eeq

\medskip

It was already clear from \eqref{primastima} that the \omegalimit\ $\omega$ is non-empty.
Indeed, the continuous $\calV$-valued function $\Rho$ is also bounded,
so that there exists a sequence $\tn\nearrow+\infty$ such that
the sequence $\graffe{\Rho(\tn)}$ is weakly convergent in~$\calV$.
More precisely, any sequence of times that tends to infinity contains a subsequence of this type.
Thus, it remains to prove the second part of the statement.
Therefore, we fix an element $\Rhoo\in\omega$ 
and a corresponding sequence $\{\tn\}$ like in the definition~\eqref{omegalim}.
We also fix $T\in(0,+\infty)$ and set, \aat,
\Bsist
  && \mun(t) := \mu(\tn+t), \quad
  \rhon(t) := \rho(\tn+t), \quad
  \zetan(t) := \zeta(\tn+t)
  \non
  \\
  && \muGn(t) := \muG(\tn+t), \quad
  \rhoGn(t) := \rhoG(\tn+t), \quad
  \zetaGn(t) := \zetaG(\tn+t)
  \non
  \\
  && \un(t) := u(\tn+t)
  \non
\Esist
and notice that \eqref{hpu} and the interpolation inequality \eqref{interpolaz} 
applied to $\un$ imply that
\Beq
  \un \to 0
  \quad \hbox{strongly in $\L\infty\Ldue$}.
  \label{convun}
\Eeq
Moreover, it is clear that the 6-tuple $\soluzn$
satisfies the regularity conditions \Regsoluz\
and the equations \accorpa{prima}{terza} with $u$ replaced by~$\un$,
as well as the initial condition $\rhon(0)=\rho(\tn)$.
In particular, by construction, we have~that
\Beq
  (\rhon,\rhoGn)(0) \to \Rhoo
  \quad \hbox{weakly in $\calV$}.
  \label{limrhonz}
\Eeq
Furthermore, the global estimates already performed on $\soluz$
immediately imply some estimates on $\soluzn$ that are uniform with respect to~$n$.
Here is a list.
From \eqref{primastima} and \eqref{secondastima}, we infer that
\Beq
  \norma\Rhon_{\H1\calH\cap\L\infty\calV} \leq c \,.
  \label{stimarhon}
\Eeq
By \eqref{globint}, we also deduce that
\Bsist
  && (\nabla\mun,\nablaG\muGn) \to 0
  \quad \hbox{strongly in $(\L2\calH)^3$}
  \label{nablamun}
  \\
  && (\dt\rhon,\dt\rhoGn) \to 0
  \quad \hbox{strongly in $\L2\calH$}.
  \label{dtrhon}
\Esist
On the other hand, \eqref{terzastima} yields a uniform estimate on the mean value $\mean\Mun$.
By combining this with \eqref{nablamun}, we conclude that
\Beq
  \norma\Mun_{\L2\calV} \leq \CT \,.
  \label{stimamun}
\Eeq
In the next steps, we perform further estimates that ensure
some \juerg{additional convergence properties} for $\soluzn$ on the interval $(0,T)$.

\step 
First auxiliary estimate

Once again, we treat $\beta$ and $\betaG$ as single-valued functions
and write $\beta(\rhon)$ and $\betaG(\rhon)$ in place of $\zetan$ and~$\zetaGn$, respectively.
We test \eqref{seconda}, written for $\soluzn$ at the time $s$, by $(\beta(\rhon),\beta(\rhoGn))(s)\in\calV$
and integrate over $(0,t)$ with any $t\in(0,T)$.
We have
\Bsist
  && \frac \tauO 2 \iO \Beta(\rhon(t))
  + \frac \tauG 2 \iG \BetaG(\rhoGn(t))
  + \intQt \beta'(\rhon) |\nabla\rhon|^2
  + \intSt \beta'(\rhon) |\nablaG\rhoGn|^2
  \non
  \\
  && \quad {}
  + \intQt |\beta(\rhon)|^2
  + \intSt \betaG(\rhoGn) \, \beta(\rhoGn)
  \non
  \\
  && = \intQt \bigl( \mun - \pi(\rhon) \bigr) \, \beta(\rhon)
  + \intSt \bigl( \muGn - \piG(\rhoGn) \bigr) \, \beta(\rhoGn)
  \non
  \\
  && \quad {}
  + \frac \tauO 2 \iO \Beta(\rho(\tn))
  + \frac \tauG 2 \iG \BetaG(\rhoG(\tn)).
  \label{perprimaaux}
\Esist
All of the integrals on the \lhs\ are nonnegative but the last one,
which we treat in the following way:
we notice that \eqref{hpCC}, and the fact that 
$\beta$ and $\betaG$ have the same sign by \eqref{hpBeta}, 
imply that\,
$\beta(r)\betaG(r)\geq(2\eta)^{-1}\beta^2(r)-c$\, for all~$r$.
It follows that
\Beq
  \intSt \betaG(\rhoGn) \, \beta(\rhoGn)
  \geq \frac 1 {2\eta} \intSt |\beta(\rhoGn)|^2  - \CT \,.
  \non
\Eeq
Let us come to the \rhs.
The first integral can obviously be dealt with the Young inequality.
We treat the second one this~way:
\begin{align*}
  &\pier{\intSt \bigl( \muGn - \piG(\rhoGn) \bigr) \, \beta(\rhoGn)}
  \\
  &\quad   \leq \frac 1 {4\eta} \intSt |\beta(\rhoGn)|^2
  + c \intSt |\muGn|^2
  \pier{{}+ c \intSt \bigl(1 +|\rhoGn|^2\bigr)}
   \leq \frac 1 {4\eta} \intSt |\beta(\rhoGn)|^2
  + c \,,
  \non
\end{align*}
thanks to \eqref{stimamun}.
Finally, the last two term of \eqref{perprimaaux} are bounded by \eqref{stimeBeta}.
By combining all these inequalities, we derive that
\Beq
  \norma\zetan_{\pier{\L2H}} \leq \CT\,, 
  \label{primaaux}
\Eeq
as well as an estimate for $\norma{\beta(\rhoGn)}_{\L2\HG}$ as a by-product.

\step
Second auxiliary estimate

We apply \cite[Lem.~3.1]{CGS13} \Aet\ to \eqref{seconda}, 
written for $\soluzn$ in the form
\Bsist
  && \iO \nabla\rhon \cdot \nabla v
  + \iG \nablaG\rhoGn \cdot \nablaG\vG
  + \iG \betaG(\rhoGn) \vG
  \non
  \\
  && = \iO \mun v 
  + \iG \muGn \vG  
  - \tauO \iO \dt\rhon \, v
  - \tauG \iG \dt\rhoGn \, \vG
  - \iO \bigl( \zetan + \pi(\rhon) \bigr) v
  - \iG \piG(\rhoGn) \vG .
  \non
\Esist
We obtain, in particular, that
\Beq
  \norma{\betaG(\rhoGn(t))}\juerg{_{H_\Gamma}}
  \leq c \, \bigl(
    \normaHH{\Mun(t)}
    + \pier{\normaHH{\dt\Rhon(t)}}
    + \normaH{(\zetan+\pi(\rhon))(t)}
  \bigr)
  \non
\Eeq
\Aat, where $c$ depends only on~$\Omega$.
By accounting for \eqref{stimarhon}, \eqref{stimamun} and \eqref{primaaux}, 
we conclude that
\Beq
  \norma\zetaGn_{\L2\HG} \leq \CT \,.
  \label{secondaaux}
\Eeq

\step
Limits

We collect the estimates \eqref{stimarhon}, \eqref{stimamun}, \eqref{primaaux} and \eqref{secondaaux}
and use standard weak and weak star compactness results.
For a subsequence, still labeled by~$n$, we~have
\Bsist
  & \Rhon \to \Rhoi
  & \quad \hbox{weakly star in $\H1\calH\cap\L\infty\calV$},
  \label{convRhon}
  \\
  & \Mun \to \Mui
  & \quad \hbox{weakly in $\L2\calV$},
  \label{convMun}
  \\
  & \Zetan \to \Zetai
  & \quad \hbox{weakly in $\L2\calH$}.
  \label{convZetan}
\Esist
We now prove that $\soluzi$ satisfies 
the integrated version of \accorpa{prima}{seconda}, where we read $u=0$,
with time-dependent test functions $(v,\vG)\in\L2\calV$,
and that \eqref{terza} holds true as well.
First of all, we notice that $\rhon$ converges to $\rhoi$ 
weakly star in $\L\infty{\Lx6}$, 
by \eqref{convRhon} and the continuous embedding $V\subset\Lx6$.
\gianni{Owing to~\eqref{utozero}, we see that $\un$ tends to zero strongly in $\pier{(}\L\infty{\Lx{3/2}}\pier{)^3}$
and deduce that $\rhon\un$ converges to zero weakly star in $\pier{(}\L\infty{\Lx{6/5}}\pier{)^3}$}.
Next, from \eqref{convRhon} we derive the strong convergence
\Beq
  \Rhon \to \Rhoi
  \quad \hbox{strongly in $\C0\calH$},
  \label{strongRhon}
\Eeq
\juerg{which follows from the compact embedding $\calV\subset\calH$ 
and from} applying, e.g., \cite[Sect.~8, Cor.~4]{Simon}.
Then, $(\pi(\rhon),\piG(\rhoGn))$ converges to $(\pi(\rhoi),\piG(\rhoGi))$
strongly in the same space, by \Lip\ continuity.
This concludes the proofs that \eqref{prima} and \eqref{seconda}
hold for the limiting 6-tuple in an integrated form,
which is equivalent to the pointwise formulation.
In order to derive~\eqref{terza}, i.e.,
$\zetai\in\beta(\rhoi)$ and $\zetaGi\in\betaG(\rhoGi)$,
we combine the weak convergence \eqref{convZetan} 
with the strong convergence \eqref{strongRhon}
and apply, e.g., \cite[Lemma~2.3, p.~38]{Barbu}.

\step
Conclusion

It remains to prove that the above limit leads to a stationary solution 
with the properties specified in the statement.
To this end, we first derive that $\Rhoi$ belongs to $\L2\calW$ 
and solves the boundary value problem
\Bsist
  && - \Delta\rhoi + \zetai + \pi(\rhoi) = \mui
  \quad \aeQ\,,
  \label{bvpO}
  \\
  && \dn\rhoi - \DeltaG\rhoGi + \zetaGi + \piG(\rhoGi) = \muGi
  \quad \aeS \,.
  \label{bvpG}
\Esist
From \eqref{dtrhon} \pier{and \eqref{convRhon}}, we see that $\dt\Rhoi$ vanishes identically.
Thus, we are dealing with a time-dependent elliptic problem
in a variational form
and can use the following well-known estimates from trace theory
and from the theory of elliptic equations.
For any $v$ and $\vG$ that make the \rhs s meaningful,
we have that
\Bsist
  && \norma{\dn v}_{\HxG{-1/2}}
  \leq \CO \bigl( \norma v_{\Huno} + \norma{\Delta v}_{\Ldue} \bigr)\,,
  \label{dnHmum}
  \\[1mm]
  && \norma{\dn v}_{\LdueG}
  \leq \CO \bigl( \norma v_{\Hx{3/2}} + \norma{\Delta v}_{\Ldue} \bigr)\,,
  \label{dnLd}
  \\[1mm]
  && \norma v_{\Hdue} \leq \CO \bigl( \norma{v\suG}_{\HxG{3/2}} + \norma{\Delta v}_{\Ldue} \bigr)\,,
  \label{stimaHdueD}
  \\[1mm]
  && \norma\vG_{\HdueG} \leq \CO \bigl( \norma\vG_{\HunoG} + \norma{\DeltaG\vG}_{\LdueG} \bigr)\,,
  \label{stimaHdueG}
  \\[1mm]
  && \norma\vG_{\HxG{3/2}} \leq \CO \bigl( \norma\vG_{\HunoG} + \norma{\DeltaG\vG}_{\HxG{-1/2}} \bigr)\,,
  \label{stimaHtmG}
\Esist
where the positive constant $\CO$ depends only on~$\Omega$.
By taking test functions $(v,0)$ with $v\in\Hunoz$,
we derive that \eqref{bvpO} holds in the sense of distributions on~$\QT$.
This implies that $\Delta\rhoi\in\L2H$,
so that $\dn\rhoi$ is a well-defined element of $\L2{\HxG{-1/2}}$ (see~\eqref{dnHmum})
satisfying the integration--by--parts formula in a generalized sense.
Coming back to \eqref{seconda} written for our limiting solution 
and arbitrary test functions $(v,\vG)\in\calV$, 
we deduce that \eqref{bvpG} holds in a generalized sense.
From this, we infer that $\DeltaG\rhoGi\in\L2{\HxG{-1/2}}$,
so that $\rhoGi\in\L2{\HxG{3/2}}$ (see~\eqref{stimaHtmG}).
It follows that $\rhoi\in\L2W$ by~\eqref{stimaHdueD}.
In particular, $\dn\rhoi\in\L2\HG$ by \eqref{dnLd} so that $\DeltaG\rhoGi\in\L2\HG$
and $\rhoGi\in\L2\WG$ by~\eqref{stimaHdueG}.
Finally, as all of the ingredients are $L^2$ functions,
it is also clear that equations \accorpa{bvpO}{bvpG} hold almost everywhere.

At this point, we are ready to conclude.
Since $\dt\Rhoi$ vanishes and the same holds for $(\nabla\mui,\nablaG\muGi)$ by~\eqref{nablamun}, 
there exist $\Rhos\in\calV$ and $\musi\in L^2(0,T)$ such that
\Beq
  \Rhoi(x,t) = \Rhos(x)
  \aand
  \Mui(x,t) = (\musi(t),\musi(t))
  \quad \hbox{for a.a.\ $(x,t)\in\QT$} \,.
  \non
\Eeq
We prove that $\Zetai$ is time independent as well and that $\musi$ is a constant,
by accounting for our assumptions on the graphs $\beta$ and~$\betaG$:
at least one of them is single-valued.
Assume first that $\beta$ is single-valued.
Then $\zetai=\beta(\rhoi)$ is time independent and attains
the value $\zetas:=\beta(\rhos)$ at any time.
From~\eqref{bvpO}, it follows that $\mui$ is time independent as well,
so that the function $\musi$ is a constant that we term~$\mus$.
Thus, the \rhs\ of \eqref{bvpG} is the same constant~$\mus$.
As this does not depend on time, the same holds for~$\pier{\zetaGi}$.
Hence, it attains some value $\zetaGs\in\HG$ \aet.
Assume now that $\betaG$ is single-valued.
Then, we first use \eqref{bvpG} to derive that 
$\zetaGi=\betaG(\rhoGi)$ and $\muGi$ are time independent.
In particular, $\pier{\mui}$~takes some constant value~$\mus$,
so that $\zetai$ is time independent, by comparison in~\eqref{bvpO}.
As in both cases \eqref{terzas} holds as a consequence of \eqref{terza} for the limiting solution,
the quadruplet $\soluzs$ is a stationary solution
corresponding to the value $\mus$ of the chemical potential.
Finally, we prove that $\Rhos$ coincides with the given~$\Rhoo$.
Indeed, \eqref{convRhon} implies weak convergence also in $\C0\calH$, whence
\Beq
  \Rhon(0) \to \Rhoi(0) = \Rhos
  \quad \hbox{weakly in~$\calH$}.
  \non
\Eeq
By comparison with \eqref{limrhonz}, we conclude that $\Rhos=\Rhoo$,
and the proof is complete.


\section*{Acknowledgments}
PC and GG gratefully acknowledge some financial support \pier{from 
the MIUR-PRIN Grant 2015PA5MP7 ``Calculus of Variations'',}
the GNAMPA (Gruppo Nazionale per l'Analisi Matematica, 
la Probabilit\`a e le loro Applicazioni) of INdAM (Isti\-tuto 
Nazionale di Alta Matematica) and the IMATI -- C.N.R. Pavia.


\vspace{3truemm}


\vspace{3truemm}

\Begin{thebibliography}{10}

\pier{%
\bibitem{bai}
F.~Bai, C.M. Elliott, A.~Gardiner, A.~Spence and A.M. Stuart,  
The viscous {C}ahn-{H}illiard equation. {I}. {C}omputations, 
{\it Nonlinearity\/}  {\bf 8} (1995) 131-160.}

\bibitem{Barbu}
V. Barbu,
``Nonlinear Differential Equations of Monotone Type in Banach Spaces'',
Springer,
London, New York, 2010.

\bibitem{BL}
J. Bergh and J. L\"ofstr\"om,
``Interpolation spaces. An introduction'',
Grundlehren der mathematischen Wissenschaften {\bf 223},
Springer-Verlag, Berlin-New York, 1976.

\pier{\bibitem{CahH} 
J.W. Cahn and J.E. Hilliard, 
Free energy of a nonuniform system I. Interfacial free energy, 
{\it J. Chem. Phys.\/}
{\bf 2} (1958) 258-267.}

\bibitem{CaCo}
L. Calatroni and P. Colli,
Global solution to the Allen--Cahn equation with singular potentials and dynamic boundary conditions,
{\it Nonlinear Anal.\/} {\bf 79} (2013) 12-27.

\pier{%
\bibitem{CGM13}
 L.\ {C}herfils, S.\ {G}atti and A.\ {M}iranville, 
\newblock A variational approach to a {C}ahn--{H}illiard model in a domain 
with nonpermeable walls,
\newblock J.\ Math.\ Sci.\ (N.Y.) {\bf 189} (2013) 604-636. 
\bibitem{CMZ11}
 L.\ {C}herfils, A.\ {M}iranville and S.\ {Z}elik,
The {C}ahn--{H}illiard equation with logarithmic potentials,
{\it Milan J. Math.\/} {\bf 79} (2011) 561-596.%
\bibitem{CP14}
L.\ {C}herfils and M.\ {P}etcu, 
	\newblock A numerical analysis of the {C}ahn--{H}illiard equation with non-permeable walls,
	\newblock Numer.\ Math. {\bf 128} (2014) 518-549. 
\bibitem{CFP} 
R. Chill, E. Fa\v sangov\'a and J. Pr\"uss,
Convergence to steady states of solutions of the Cahn-Hilliard equation with dynamic boundary conditions,
{\it Math. Nachr.\/} 
{\bf 279} (2006) 1448-1462.
\bibitem{CF1} 
P.\ {C}olli and T.\ {F}ukao, 
{C}ahn--{H}illiard equation with dynamic boundary conditions 
and mass constraint on the boundary, 
{\it J. Math. Anal. Appl.\/} {\bf 429} (2015) 1190-1213.
\bibitem{CF2} 
P.\ {C}olli and T.\ {F}ukao, 
Equation and dynamic boundary condition of 
Cahn--Hilliard type with singular potentials, 
{\it Nonlinear Anal.\/} {\bf 127} (2015) 413-433.
\bibitem{CGPS3} 
P. Colli, G. Gilardi, P. Podio-Guidugli and J. Sprekels,
Well-posedness and long-time behaviour for 
a nonstandard viscous Cahn-Hilliard system, 
{\it SIAM J. Appl. Math.} {\bf 71} (2011) 1849-1870.
\bibitem{CGS3}
P. Colli, G. Gilardi and J. Sprekels,
On the Cahn--Hilliard equation with dynamic 
boundary conditions and a dominating boundary potential,
{\it J. Math. Anal. Appl.\/} {\bf 419} (2014) 972-994.
\bibitem{CGS5}
P. Colli, G. Gilardi and J. Sprekels,
A boundary control problem for the pure Cahn-Hilliard equation
with dynamic boundary conditions,
{\it Adv. Nonlinear Anal.\/} {\bf 4} (2015) 311-325.
\bibitem{CGS4}
P. Colli, G. Gilardi and J. Sprekels,
A boundary control problem for the viscous Cahn-Hilliard equation
with dynamic boundary conditions,
{\it Appl. Math. Optim.\/} {\bf 73} (2016) 195-225.
}

\bibitem{CGS13}
P. Colli, G. Gilardi, J. Sprekels,
On a Cahn--Hilliard system with convection and
dynamic boundary conditions,  \pier{{\it Ann. Mat. Pura Appl. (4),} DOI~10.1007/s10231-018-0732-1
(see also preprint arXiv:1704.05337 [math.AP] (2017), pp.~1-34).}%

\pier{\bibitem{CGS14}
P. Colli, G. Gilardi, J. Sprekels,
Optimal velocity control of a viscous Cahn--Hilliard system with convection and dynamic boundary conditions, 
{\it SIAM J. Control Optim.\/}, to appear (2018)
(see also preprint arXiv:1709.02335 [math.AP] (2017), pp.~1-28).
\pier{\bibitem{CGS15}
P. Colli, G. Gilardi, J. Sprekels,
Optimal velocity control of a convective Cahn--Hilliard system with double obstacles
and dynamic boundary conditions: a `deep quench' approach, 
preprint arXiv:1709.03892 [math.AP] (2017), pp.~1-30.}
\bibitem{CS}
P. Colli and J. Sprekels,
Optimal control of an Allen--Cahn equation 
with singular potentials and dynamic boundary condition,
{\it SIAM J. Control Optim.\/} {\bf 53} (2015) 213-234.%
\bibitem{EllSt} 
C.M. Elliott and A.M. Stuart, 
Viscous Cahn--Hilliard equation. II. Analysis, 
{\it J. Differential Equations\/} 
{\bf 128} (1996) 387-414.
\bibitem{EllSh} 
C.M. Elliott and S. Zheng, 
On the Cahn--Hilliard equation, 
{\it Arch. Rational Mech. Anal.\/} 
{\bf 96} (1986) 339-357.
\bibitem{FG} 
E. Fried and M.E. Gurtin, 
Continuum theory of thermally induced phase transitions based on an order 
parameter, {\it Phys. D} {\bf 68} (1993) 326-343.
\bibitem{FY} 
T. Fukao and N. Yamazaki, 
A boundary control problem for the equation 
and dynamic boundary condition of Cahn--Hilliard type, 
{in ``Solvability, 
Regularity, Optimal Control of Boundary Value Problems for PDEs'', 
P.~Colli, A.~Favini, E.~Rocca, G.~Schimperna, J.~Sprekels~(ed.), 
Springer INdAM Series {\bf 22}, Springer, Milan, 2017, pp.~255-280.}
\bibitem{G1} 
C.G. Gal, A Cahn-Hilliard model in bounded domains with permeable walls. 
{\it Math. Methods Appl. Sci.\/} {\bf 29} (2006) 2009-2036.
\bibitem{GW} 
C.G. Gal and H. Wu, Asymptotic behavior of a Cahn-Hilliard equation with Wentzell boundary conditions and mass conservation. {\it Discrete Contin. Dyn. Syst.\/} {\bf 22} (2008) 1041-1063.%
}

\bibitem{GiMiSchi} 
G. Gilardi, A. Miranville and G. Schimperna,
On the Cahn--Hilliard equation with irregular potentials and dynamic boundary conditions,
{\it Commun. Pure Appl. Anal.\/} 
{\bf 8} (2009) 881-912.

\pier{%
\bibitem{GiMiSchi2} 
G. Gilardi, A. Miranville and G. Schimperna,
Long-time behavior of the Cahn-Hilliard 
equation with irregular potentials 
and dynamic boundary conditions,
{\it Chin. Ann. Math. Ser. B\/} 
{\bf 31} (2010) 679-712.%
}

\bibitem{GiRo}
G. Gilardi, E. Rocca,
Well posedness and long time behaviour for a singular phase field system of conserved type,
{\em IMA J. Appl. Math.} {\bf 72} (2007) 498-530.

\pier{%
\bibitem{GMS11} 
	 G.R.\ {G}oldstein, A.\ {M}iranville and G.\ {S}chimperna, 
	\newblock A {C}ahn--{H}illiard model in a domain with non-permeable walls,
	\newblock Phys. D {\bf 240} (2011) 754-766.
\bibitem{GM13} 
	 G.R.\ {G}oldstein and A.\ {M}iranville,
	\newblock A {C}ahn--{H}illiard--{G}urtin model with dynamic boundary conditions,
	\newblock Discrete Contin.\ Dyn.\ Syst.\ Ser.\ S {\bf 6} (2013) 387-400.
\bibitem{Gu} 
M. Gurtin, 
Generalized Ginzburg-Landau and Cahn-Hilliard equations based on a microforce balance,
{\it Phys.~D\/} {\bf 92} (1996) 178-192.
\bibitem{Kub12} 
	 M.\ {K}ubo,
	\newblock The {C}ahn--{H}illiard equation with time-dependent constraint, 
	\newblock Nonlinear Anal. {\bf 75} (2012) 5672-5685. 
\pier{\bibitem{LW}
C. Liu and H. Wu,
An energetic variational approach for the Cahn--Hilliard equation with dynamic boundary conditions: derivation and analysis,  preprint arXiv:1710.08318 [math.AP] (2017), pp.~1-68.}
\bibitem{MZ}
A. Miranville and S. Zelik,
Robust exponential attractors for {C}ahn-{H}illiard type
equations with singular potentials,
{\it Math. Methods Appl. Sci.} {\bf 27} (2004) 545-582.
\bibitem{NovCoh}
A.~Novick-Cohen, On the viscous {C}ahn-{H}illiard equation, in
``Material instabilities in continuum mechanics'' ({E}dinburgh, 1985--1986),
Oxford Sci. Publ., Oxford Univ. Press, New York, 1988, pp.~329-342.
\bibitem{Podio}
P. Podio-Guidugli, 
Models of phase segregation and diffusion of atomic species on a lattice,
{\it Ric. Mat.} {\bf 55} (2006) 105-118.
\bibitem{PRZ} 
J. Pr\"uss, R. Racke and S. Zheng, 
Maximal regularity and asymptotic behavior of solutions for the 
Cahn--Hilliard equation with dynamic boundary conditions,  
{\it Ann. Mat. Pura Appl.~(4)\/}
{\bf 185} (2006) 627-648.
\bibitem{RZ} 
R. Racke and S. Zheng, 
The Cahn--Hilliard equation with dynamic boundary conditions, 
{\it Adv. Differential Equations\/} 
{\bf 8} (2003) 83-110.%
}
\bibitem{Simon}
J. Simon,
Compact sets in the space $L^p(0,T; B)$,
{\it Ann. Mat. Pura Appl.~(4)\/} 
{\bf 146} (1987) 65-96.

\pier{%
\bibitem{WZ} H. Wu and S. Zheng,
Convergence to equilibrium for the Cahn--Hilliard equation 
with dynamic boundary conditions, {\it J. Differential Equations\/}
{\bf 204} (2004) 511-531.%
}

\End{thebibliography}

\End{document}


  \\
  \separa
  && \iO \nabla w \cdot \nabla\xi
  + \iG \nablaG\wG \cdot \nablaG\xiG
  = \frac 12 \, \frac d{dt} \, \normaHH{(w,\wG)}^2
  \non
  \\
  && \quad \hbox{if $(w,\wG)\in\L2\calV$, \ $\dt(w,\wG) \in \L2\calVsz$}
  \non
  \\   
  && \aand
  \hbox{$\Xi = \calN(\dt(w,\wG))$} \,.
  && \< \gstar , \calN\gstar >_{\calV}
  = \norma\gstar_*^2
  \qquad \hbox{if $\gstar \in \calVsz$},
  \\
  \separa
  && \iO \nabla w \cdot \nabla\xi
  + \iG \nablaG\wG \cdot \nablaG\xiG
  = \normaHH{(w,\wG)}^2
  \non
  \\
  && \quad \hbox{if $(w,\wG) \in \calHz$ and $(\xi,\xiG) = \calN(w,\wG)$} \,.
  \\
  \separa